\magnification 1200
\baselineskip 14pt
 1
\font\large=cmbx12 scaled \magstep 1

\vsize 19.5cm
\hsize 14cm
\def\qed{\hfill\vbox{\hrule\hbox{\vrule\kern3pt
\vbox{\kern6pt}\kern3pt\vrule}\hrule}\bigskip}
\input amssym

\def\normalization{1}
\def\technical{2}
\def\existencemaximals{3}
\def\finiteextension{4}
\def\classicalweak{5}
\def\intersectionmaximals{6}
\def\classicalstrong{7}

\centerline{\large The Nullstellensatz without the Axiom of Choice}
\smallskip
\centerline{Enrique Arrondo}

\bigskip

\noindent{\bf Abstract.} We give a short proof of the most general version of the Nullstellensatz without using the Axiom of Choice.

\bigskip


\bigskip

\noindent The general theory of schemes makes a strong use of the main two results of Commutative Algebra that require the Axiom of Choice: the existence of maximal ideals and the fact that the radical of an ideal is the intersection of all prime ideals containing the given ideal. However, when one deals with finitely generated algebras over a field $K$, these two results are a consequence of the Nullstellensatz (moreover, in this case, the radical of an ideal is the intersection of the maximal ideals containing it; hence finitely generated $K$-algebras are what is called Jacobian rings). In fact, they can be considered respectively as extensions of the so-called weak and strong versions of the Nullstellensatz, when the ground field is an arbitrary field (i.e. not necessarily algebraically closed, which is the assumption in the classical version of the Nullstellensatz).

It is thus natural to look for proofs of the Nullstellensatz (in any possible version) avoiding any use of the Axiom of Choice. This has been approached in a more restrictive way: finding constructive proofs (see [2] Corollary VI-1.12 or [3] Theorem V1-3.5 for the most general version we are going to deal with). However, these proofs still keep the traditional approach of using the notion of integral extensions of rings.

The goal of this note is to find simple and short proofs of all these versions of the Nullstellensatz, without any use of the Axiom of Choice. The main idea is taken from [1], in which we gave an elementary proof of the (classical) Weak Nullstellensatz for algebraically closed fields. We state that main idea separately in Lemma \technical, and it will be the keystone for our main proofs. We want to stress that --forgetting about using or not the Axiom of Choice-- we do not know of any previous simple proof of the most general version of the Nullstellensatz.

In the case of an algebraically closed field $K$, the Nullstellensatz implies that maximal ideals of a finitely generated $K$-algebra are in bijection with the closed points of the corresponding affine variety. When $K$ is arbitrary, the generalized Nullstellensatz implies that maximal ideals correspond now to sets of conjugate points in finite extensions of $K$ (this is the spirit of Theorem \finiteextension). Hence, a generalized Weak Nullstellensatz should state the existence of maximal ideals, while a generalized Strong Nullstellensatz should state that the radical of an ideal is the intersection of the maximal ideals containing the ideal. We will prove all these results without using the Axiom of Choice. Although some of the proofs are well-known, we will include them here for the sake of completeness. 

We want to thank Henri Lombardi for indicating us the references [2] and [3]. This research was developed in the framework of the project MTM2015-65968-P funded by the Spanish Government.

\bigskip

We start with the generalized Noether Normalization Lemma. We reproduce here the standard proof when $K$ is arbitrary (the case when $K$ is infinite is simple and it is also standard, and can be found in [1]).

\proclaim Lemma \normalization. If $K$ be any field and $f$ is a
nonconstant polynomial in $K[x_1,\ldots,x_n]$ with $n\ge2$, then there is an automorphism of $K[x_1,\ldots,x_n]$ transforming $f$ into a monic polynomial in the varible $x_n$.

\noindent{\it Proof.} Let $m$ be larger than any exponent $i_j$ of any monomial $x_1^{i_1}\dots x_n^{i_n}$ appearing in $f$. An automorphism of $K[x_1,\ldots,x_n]$ leaving invariant $x_n$ and mapping each other $x_i$ to $x_i+x_n^{m^{n-1-i}}$ sends any monomial  $x_1^{i_1}\dots x_n^{i_n}$ to a polynomial whose monomial of highest degree is $x_n^{i_1m^{n-1}+\dots+i_{n-1}m+i_n}$. Regarding this degree as a number written in base $m$, it follows that $i_1m^{n-1}+\dots+i_{n-1}m+i_n\le i'_1m^{n-1}+\dots+i'_{n-1}m+i'_n$ if and only if $(i_1,\dots,i_{n-1},i_n)\le(i'_1,\dots,i'_{n-1},i'_n)$ in the lexicographical order. Therefore, such automorphism maps $f$ to a monic polynomial in $x_n$ of degree $i_1m^{n-1}+\dots+i_{n-1}m+i_n$, where $(i_1,\dots,i_{n-1},i_n)$ is the maximum of the set of exponents of $f$ when ordered lexicographically.
\qed

The main technical lemma we will need is the following (which keeps the main idea of [1], and we essentially copy the proof there adapted to our more general context).

\proclaim Lemma \technical. Let $I\subset K[x_1,\dots,x_n]$ be a proper ideal containing a polynomial $g$ that is monic in the variable $x_n$. If $M'\subset K[x_1,\dots,x_{n-1}]$ is a proper ideal containing $I\cap K[x_1,\dots,x_{n-1}]$, then the ideal of $K[x_1,\dots,x_n]$ generated by $I$ and $M'$ is proper.

\noindent{\it Proof.} Suppose for contradiction that there exist $f$ in $I$ and $f'$ generated by $M'$ such $1=f+f'$. Thus we can write
$f=f_0+f_1x_n+\ldots+f_dx_n^d$, with all the $f_i$ in
$K[x_1,\ldots,x_{n-1}]$, and such that $f_0-1,f_1,\dots,f_d\in M'$. On the other hand, we can express the
monic polynomial $g$ in the form
$g=g_0+g_1x_n+\ldots+g_{e-1}x_n^{e-1}+x_n^e$ with $g_j$ in
$K[x_1,\ldots,x_{n-1}]$ for $j=0,\ldots,e-1$.

Let $R$ be the resultant of $f$ and
$g$ with respect to the variable $x_n$. In other words, $R$ is
the polynomial in $K[x_1,\ldots,x_{n-1}]$ given by the determinant

$$R=\left|\matrix{
f_0&f_1&\ldots&f_d&0&0&\ldots&0\cr
0&f_0&\ldots&f_{d-1}&f_d&0&\ldots&0\cr
&&\ddots\cr
0&\ldots&0&f_0&f_1&\ldots&f_{d-1}&f_d\cr
g_0&g_1&\ldots&g_{e-1}&1&0&\ldots&0\cr
0&g_0&\ldots&g_{e-2}&g_{e-1}&1&0\ldots&0\cr
&&\ddots&&&&\ddots\cr
0&\ldots&0&g_0&g_1&\ldots&g_{e-1}&1}
\right| 
\matrix{\left.\matrix{\cr\vphantom{\ddots}\cr\cr\cr}\right\}
e\
$rows$\cr
\left.\matrix{\cr\vphantom{\ddots}\cr\cr\cr}\right\}d\
$rows$}\hskip -1.1cm\matrix{\cr\cr\cr\cr\cr\cr\cr\cr.}$$
We recall the well known result that $R$ belongs to $I$: in the above determinant defining $R$, add to the first column the second one multiplied by $x_n$, then the third column multiplied by $x_n^2$, and so on until one adds the last column multiplied by $x_n^{d+e-1}$; developing the resulting determinant by the first column shows that $R$ is a linear combination of $f$ and $g$, so it is in $I$. Therefore $R$ is a member of $I\cap K[x_1,\dots,x_{n-1}]$, and hence of $M'$. But a direct inspection of the determinant defining the resultant shows that, when quotienting modulo $M'$, it reduces to the determinant of a lower-triangular matrix whose entries on its main diagonal are all $1$s. Hence $R$ is not in $M'$, which is a contradiction. \qed

From this, we can obtain the main versions of the Nullstellenstaz:

\proclaim Theorem \existencemaximals. Any proper ideal  $I\subset K[x_1,\ldots,x_n]$ is contained in a maximal ideal.

\noindent{\it Proof.} Let us assume $I\neq0$, since
otherwise the result is trivial. We prove the theorem
by induction on $n$. The case $n=1$ is immediate, because any
nonzero proper ideal $I$ of $K[x]$ is generated by a
nonconstant polynomial. The ideal generated by any irreducible factor of such polynomial is a maximal ideal containing $I$.

We assume now $n>1$. Lemma 1 allows us to suppose that $I$ contains a polynomial $g$ monic in the variable $x_n$. Fixing such a polynomial $g$, we consider the ideal $I':=I\cap K[x_1,\ldots,x_{n-1}]$. Since
$1$ is not in $I$, it follows that $I'$ is a proper ideal.
Therefore, by the induction hypothesis there is a maximal ideal $M'$ of $K[x_1,\ldots,x_{n-1}]$ containing $I'$. 

If we consider the field $K'=K[x_1,\ldots,x_{n-1}]/M'$, there exists a natural surjective map $K'[x_n]\to K[x_1,\ldots,x_{n}]/\big(I+(M')\big)$. By Lemma \technical, its kernel $J$ is a proper ideal of $K'[x_n]$, hence $K'[x_n]/J$, and thus also $K[x_1,\ldots,x_{n}]/\big(I+(M')\big)$ possesses a maximal ideal. The pullback to $K[x_1,\dots,x_n]$ of such maximal ideal yields a maximal ideal containing $I$, which completes the proof.
\qed

\proclaim Theorem \finiteextension. A prime ideal $I\subset K[x_1,\dots,x_n]$ is maximal if and only if the quotient $K[x_1,\dots,x_n]/I$ is a finite extension of $K$.

\noindent{\it Proof.} We first prove that the finite dimension of $K[x_1,\dots,x_n]/I$ implies that $I$ is maximal. Equivalently, we need to prove that the quotient $K[x_1,\dots,x_n]/I$ is a field. To that purpose, we fix a non-zero element of the quotient, i.e. the class of a polynomial $f$ not in $I$. Since the classes of $1,f,f^2,...$ are not $K$-linearly independent modulo $I$, there is a nontrivial combination $\lambda_0+\lambda_1f+\dots+\lambda_df^d$ (with $\lambda_0,\dots,\lambda_d\in K)$ that is in $I$. Since $f$ is not in $I$ and $I$ is prime, we can assume $\lambda_0\ne0$. Hence the class of $\lambda_0^{-1}(-\lambda_1-\dots-\lambda_df^{d-1})$ is an inverse of the class of $f$. This proves that $K[x_1,\dots,x_n]/I$ is a field, and hence $I$ is a maximal ideal.

We prove now the converse by induction on $n$, the case $n=1$ being trivial, since the quotient of $K[x]$ by the ideal generated by an irreducible polynomial has finite dimension. So assume $n\ge2$ and assume also that $I$ is a maximal ideal. By Lemma \normalization, we can suppose $I$ to contain a polynomial $g$ that is monic in $x_n$. Since $I\cap K[x_1,\dots,x_{n-1}]$ is proper, because it does not contain $1$, it is contained in a maximal ideal $M'\subset K[x_1,\dots,x_{n-1}]$  (see Theorem \existencemaximals). By Lemma \technical, the ideal generated by $I$ and $M'$ is proper. Since $I$ is maximal, this implies that $M'$ is contained in $I$. In other words, $I\cap K[x_1,\dots,x_{n-1}]=M'$.

By induction hypothesis, the extension $K\subset K[x_1,\dots,x_{n-1}]/M'$ is finite. On the other hand, the extension $K[x_1,\dots,x_{n-1}]/M'\subset K[x_1,\dots,x_n]/I$ is also finite, because of the existence of a monic polynomial $g$ in $I$. This proves that $K[x_1,\dots,x_n]/I$ is finite over $K$, as wanted.
\qed

\noindent{\bf Remark \classicalweak.} When $K$ is algebraically closed, any finite extension of it is isomorphic to $K$. Hence any maximal ideal is the kernel of a morphism of $K$-algebras $K[x_1,\dots,x_n]\to K$. This morphisms is determined by the images $a_1,\dots,a_n$ of $x_1,\dots,x_n$, hence it is the evaluation at the point $(a_1,\dots,a_n)\in{\Bbb A}^n_K$. In other words, the maximal ideals of $K[x_1,\dots,x_n]$ are, in this case, the ideals of points. In this context, Theorem \existencemaximals\ is saying that, if $I\subset K[x_1,\dots,x_n]$ is a proper ideal, there is a point in ${\Bbb A}^n_K$ such that all polynomials of $I$ vanish at it. This is the classical statement of the Weak Nullstellensatz.

\noindent When $K$ is an arbitrary field, what we get is that maximal ideals are now kernels of a morphism of $K$-algebras $K[x_1,\dots,x_n]\to K'$, where $K'$ is a finite extension of $K$, i.e. they correspond to points whose coordinates are in a finite extension of $K$. Of course, any $K$-automorphism of $K'$ (i.e. any element of the Galois group of the extension) provides another morphism with the same kernel. Therefore a maximal ideal of $K[x_1,\dots,x_n]$ can be regarded as the ideal of a set of conjugate points in a finite extension of $K$.

\proclaim Theorem \intersectionmaximals. If $I\subset K[x_1,\dots,x_n]$ is a proper ideal, its radical is the intersection of all maximal ideals containing $I$.

\noindent{\it Proof.} It is clear that the radical of $I$ is contained in any maximal ideal containing $I$, so that we only need to check that a polynomial $f$ in all maximal ideals containing $I$ is in the radical of $I$. Assume, for contradiction, that $f$ is not in $\sqrt{I}$. This is equivalent to say that the $K$-algebra $\big(K[x_1,\dots,x_n]/I\big)_{\bar f}$ is not zero (where $\bar f$ stands for the class of $f$ modulo $I$). Since this $K$-algebra is isomorphic to the quotient of $K[x_1,\dots,x_n,x_{n+1}]$ by the ideal generated by $I$ and $x_{n+1}f-1$, Theorem \existencemaximals\ implies that it possesses a maximal ideal $\tilde M$. 

Consider now the kernel $P$ of the natural map $K[x_1,\dots,x_n]\to \big(K[x_1,\dots,x_n]/I\big)_{\bar f}/\tilde M$, which is a prime ideal containing $I$. Since the image of $f$ is the class of a unit, $P$ does not contain $f$. Hence, we will find the wanted contradiction if we prove that $P$ is in fact a maximal ideal. And this holds by Theorem \finiteextension, because $K[x_1,\dots,x_n]/P$ can be regarded as a subspace of $\big(K[x_1,\dots,x_n]/I\big)_{\bar f}/\tilde M$, which is finite dimensional.
\qed

\noindent{\bf Remark \classicalstrong.} According to Remark \classicalweak, Theorem \intersectionmaximals\ can be read, when $K$ is algebraically closed, as saying that the radical of an ideal $I\subset K[x_1,\dots,x_n]$ is the ideal formed by all polynomials vanishing at $V(I)$ (where $V(I)$ is the set of points of ${\Bbb A}^n_K$ that kill all polynomials of $I$). This is the classical statement of the Strong Nullstellensatz. The proof of Theorem \intersectionmaximals\ (which is the standard one) can be considered as the generalization of the classical Rabinovich trick that allows to prove the strong version of the Nullstellensatz from the weak one.

\bigskip

\noindent{\bf References.}

\medskip

\item{[1]} E. Arrondo, {\it Another Elementary Proof of the
Nullstellensatz}, Amer. Math. Monthly 113 (2006), no. 2, 169-171.

\item{[2]} H. Lombardi, C. Quitt\'e, {\it Commutative algebra: Constructive methods. Finite projective modules}, arXiv: 1605.04832.

\item{[3]} R. Mines, F. Richman, W. Ruitenburg, {\it A Course in Constructive Algebra}, Springer 1988.

\centerline{Departamento de \'Algebra, Geometr\'{\i}a y Topolog\'{\i}a}
\centerline{Facultad de Ciencias Matem\'aticas}
\centerline{Universidad Complutense de Madrid}
\centerline{28040 Madrid, Spain}
\centerline{arrondo@mat.ucm.es}

\end